\documentclass[12pt]{amsart}
\usepackage{times,fullpage}
\usepackage{latexsym,amscd,amssymb,url}
\newtheorem*{prob3}{Problem 3 {\rm (Hirzebruch-Thom)}}

\newtheorem*{prob10}{Problem 10 {\rm (Chern)}}

\newtheorem*{prob11}{Problem 11 {\rm (Chern)}}

\newtheorem*{prob12}{Problem 12 {\rm (Chern)}}

\newtheorem*{prob13}{Problem 13}

\newtheorem*{prob20}{Problem 20 {\rm (Kodaira-Spencer)}}

\newtheorem*{prob25}{Problem 25}

\newtheorem*{prob28}{Problem 28}

\newtheorem*{prob31}{Problem 31}

\theoremstyle{remark}

\theoremstyle{definition}

\newcommand{\C}{\mathbb{ C}}
\newcommand{\del}{\partial}

\newcommand{\Q}{\mathbb{ Q}}
\newcommand{\R}{\mathbb{ R}}
\newcommand{\PP}{\mathbb{ P}}
\newcommand{\HH}{\mathbb{ H}}

\newcommand{\Z}{\mathbb{ Z}}

\newcommand{\EE}{{\mathbb E}}

\newcommand{\Td}{\operatorname{Td}}

\begin{document}

\title{Updates on Hirzebruch's 1954 Problem List}
\author{D.~Kotschick}
\address{Mathematisches Institut, {\smaller LMU} M\"unchen,
Theresienstr.~39, 80333~M\"unchen, Germany}
\email{dieter@member.ams.org}

\date{May 19, 2013; \ \copyright \ D.~Kotschick 2013}
\thanks{Research done at the Institute for Advanced Study in Princeton with the support of The Fund For Math and The Oswald Veblen Fund.}

\begin{abstract}
We present updates to the problems on Hirzebruch's 1954 problem list focussing on open problems, and on those where 
substantial progress has been made in recent years. We discuss some purely topological problems, as well as geometric 
problems about (almost) complex structures, both algebraic and non-algebraic, about contact structures, and about 
(complementary pairs of) foliations.
\end{abstract}

\maketitle


\section{Introduction}

In the first week of May, 1953, a conference on fiber bundles and differential geometry was held at Cornell University.
A brief report about that conference, written by N.~E.~Steenrod, appeared in~\cite{Steenrod}, where the final two sentences 
read:

\medskip

\begin{quote}
{\it The discussions were marked by the presentations of numerous unsolved problems. 
These were recorded and a report on them is being prepared for publication.}
\end{quote}

\medskip

The young postdoc given the task of collecting the problems was F.~Hirzebruch, who was then nearing the end of the first year 
of his two-year stay at the Institute for Advanced Study (IAS) in Princeton. For him, the most important open problem was what 
he thought of as the high-dimensional Riemann-Roch problem. A few weeks after the Cornell conference, Hirzebruch proved the 
signature theorem, and, by the end of 1953, he had deduced from it what we now call the Hirzebruch--Riemann--Roch theorem.

By the end of March, 1954, when Hirzebruch submitted the collection of problems for publication, he had replaced 
the Riemann--Roch problem by several other problems that arose out of his work on the signature and Riemann--Roch theorems.
The published collection~\cite{Hir1} was reviewed in {\it Mathematical Reviews} by S.~Chern, who wrote~\cite{Ch}:

\medskip

\begin{quote}
{\it The paper lists a set of 34 unsolved problems of current interest concerned with differentiable, almost complex, and complex manifolds. 
The problems are expertly chosen; their clarification and partial or complete solutions will probably take many years and will certainly mean 
progress of the field. In stating the problems, the author tries to give their motivation, their relations to known results, and related facts. In this 
sense the paper is at the same time a resum\'e and exposition of the subject, or at least of a major part of it.}
\end{quote}

\medskip

Chern's description was certainly prescient. Reading Hirzebruch's paper~\cite{Hir1} now, 60 years after is was written, one cannot help 
noticing how modern the paper was, and how, even now, it contains nothing outdated. Hirzebruch's problems foreshadowed many of the 
developments in geometry and topology over the intervening decades. A few problems were solved right away, some were solved years 
or even decades later through major advances in the field, and a few are open to this day.

When Hirzebruch published his collected works in 1987, he compiled solutions and updates to his problems, and these appear
in~\cite[p.~762--784]{Hir2}. Some of the solutions described there draw on unpublished work of M.~Puschnigg~\cite{Puschnigg}.

This paper is an attempt to collect together further updates to Hirzebruch's problems, without repeating the information given 
in~\cite{Hir2,Puschnigg}. Thus we discuss results on Hirzebruch's problems obtained since 1987, and some earlier results that
were not mentioned in~\cite{Hir2}. We also want to advertise the problems that have remained open.
Some of the problems are quite open-ended, or admit several interpretations that are interesting. Because of this, the question of
whether a specific problem on the list is solved or not sometimes does not have a clear yes-or-no answer. Nevertheless, to 
give an overview of the material discussed below, I would like to give a few simplified pointers in this introduction, which will help
readers locate what may be of interest to them.

\subsection*{Open problems}
I believe that the following problems are essentially open:
\begin{itemize}
\item Problem~3 (elementary proof of the signature theorem),
\item Problem~10 (existence of Pfaffians of constant class),
\item Problem~13 (complex structures on $\C P^3$ and on $S^6$),
\item Problem~20 (harmonic theory for almost complex manifolds),
\item Problem~28 (complex structures on $\C P^n$).
\end{itemize}
In addition to these, the following problems, though solved to some extent, still offer interesting open questions:
\begin{itemize}
\item Problem~9 (integrability of $G$-structures),
\item Problem~11 (embeddings of $\R P^n$ in Euclidean spaces),
\item Problem~25 (complex structures on the (topological) manifold $S^2\times S^2$).
\end{itemize}

\subsection*{Progress and solutions}

I believe the only problem from~\cite{Hir1} that was wide open in 1987 and has since been solved completely, is Problem~31, 
about the topological (non-)invariance of characteristic numbers of algebraic varieties, whose solution is 
described at the end of this report. In addition, there has been important progress on the following problems:
\begin{itemize}
\item Problem~3 (geometric interpretation of the functional equation for virtual signatures),
\item Problem~10 (existence of Pfaffians of constant class),
\item Problem~11 (embeddings of $\R P^n$ in Euclidean spaces),
\item Problem~25 (complex structures on the (smooth) manifold $S^2\times S^2$).
\end{itemize}
I also found that in the case of Problems~8 and~12 some results known before 1987 were missed in~\cite{Hir2},
and so I will record them here.

\subsection*{Thanks and disclaimer}

I have known about Hirzebruch's problem list~\cite{Hir1} essentially all of my mathematical life, and I have at different times 
thought hard about several of the problems. The fact that some of them are fairly open-ended, or require interpretation,
is not a drawback, but an advantage that makes perusal of the problem list more stimulating and rewarding. 

After I gave a lecture about the solution of Problem~31 at the IAS in Princeton
in March 2013, several people asked me about the status of the other problems on this list. The present report was compiled
in part as a response to their questions. It also aims to celebrate the 60th anniversary of the list, and of Hirzebruch's
1953 work that formed the context for many of the problems. I am grateful to the IAS for its hospitality that made it possible for me
to work on Problem~31 and to compile this report. I am also grateful to D.~Davis, H.~Geiges and T.~Vogel for supplying
information on the status of certain problems.

This paper does not contain any new theorems, but rather consists of informal discussions of results, some of them my own, but
most of them due to other authors. Such informal discussions are always subjective, and are bound to be colored by the biases
and prejudices of the author. I just hope that this report does not stray too far from the taste and spirit of Hirzebruch's mathematics.
Finally, I would like to point out that I only discuss 11 of the 34 problems here, not mentioning those that were solved before
1987, unless I have something to add to Hirzebruch's report in~\cite{Hir2}.

\subsection*{Note about references}

Each section of this paper covers between one and three problems from the list from~\cite{Hir1}, and separate references
are given for each section.
Although the references given at the end of this introduction are often cited in later sections, they are not included again 
among the references for those later sections.


\vfill\eject

\newpage

\section{Two topological problems}

At the beginning of June,1953, Hirzebruch proved what he then called the ``index theorem'':
$$
\tau (M) = \langle L(p_1,\ldots,p_k),[M]\rangle \ ,
$$
where $M$ is any closed smooth oriented manifold, $\tau (M)$ is its index or signature, and $L(p_1,\ldots,p_k)$ is the $L$-polynomial
in the Pontrygain classes of $M$. Later on the term ``index'' became monopolized by the Fredholm index appearing in the 
Atiyah--Singer index theorem, and so $\tau (M)$ is now usually called the signature, and the ``index theorem'' has become the 
Hirzebruch signature theorem. 

\begin{prob3}
Give an elementary proof of the index theorem and explain the geometrical meaning of the functional equation valid for $\tau$.
\end{prob3}

Hirzebruch's proof was not elementary in the sense that it appealed to R.~Thom's determination of the oriented bordism ring over 
$\Q$, which is a rather non-elementary result in algebraic topology. Ten years later the signature theorem became a consequence 
of the Atiyah--Singer index theorem, which is even more sophisticated and ``non-elementary''. To this day no proof of the 
signature theorem is known that would qualify as ``elementary'', and so Problem~3 is wide open.

The functional equation referred to in the problem is the following. Let $\dim (M)=n$. Any homology class $x\in H_{n-2}(M;\Z)$
can be represented by a smooth closed oriented submanifold, whose signature is independent of the chosen representative.
Thus one has a well-defined virtual signature function $\tau\colon H_{n-2}(M;\Z)\longrightarrow\Z$. By taking signatures of transverse 
intersections, this can be extended to 
$$
\tau \colon H_{n-2}(M;\Z)\times\ldots \times H_{n-2}(M;\Z)\longrightarrow\Z \ .
$$
For elements $x$ and $y\in H_{n-2}(M;\Z)$ the functional equation for the virtual signature is
\begin{equation}\label{functional}
\tau (x+y)=\tau (x)+\tau (y) - \tau (x,y,x+y) \ .
\end{equation}
Hirzebruch proved this identity by using the signature theorem. A geometric interpretation of~\eqref{functional} was given by 
N.~Yu.~Netsvetaev~\cite{Nets}, who showed that if $x$, $y$ and $w=x+y$ are represented by transversely intersecting submanifolds 
$X$, $Y$ and $W$, then their bordism classes satisfy the equation
\begin{equation}\label{bordism}
[W]=[X]+[Y]-[T\tilde\times\C P^2] \ ,
\end{equation}
where the last term on the right hand side stands for a linear $\C P^2$-bundle over $T=X\cap Y\cap W$. Taking signatures and 
using the fact that $T\tilde\times\C P^2$ has the same signature as $T$, one deduces~\eqref{functional} from~\eqref{bordism}. 

Netsvetaev's proof of~\eqref{bordism} is direct and geometric, and therefore qualifies as ``elementary''. Thus his result does 
give a nice, elementary and satisfactory interpretation of~\eqref{functional}. However, this does not solve Problem~3, as 
there is still no elementary proof of the signature theorem. All that~\cite{Nets} achieved, was a geometric interpretation of the 
functional equation, and a proof of this equation that avoids the use of the signature theorem. 

\medskip

In Subsection~1.2 of~\cite{Hir1}, following Thom, Hirzebruch discusses the representability of integral homology classes in manifolds by
submanifolds. He notes that, by Thom's results, all homology classes in degrees $\leq 6$ are representable by submanifolds,
as are those of codimension at most $2$ -- this was used above to define the virtual signatures. 
Thus degree $7$ is the smallest degree in which there could be a homology class
not representable by a submanifold, and the ambient manifold would have to have dimension at least $10$. In fact, Thom had 
already constructed an example where the ambient manifold had dimension $14$. 

Problem~8 in~\cite{Hir1}, attributed to Thom, aimed at determining the smallest dimension in which there could be a manifold with an integral 
homology class not representable by a  submanifold. However, the actual formulation of that problem was more technical, 
focussing on a particular integral Steenrod power operation instead, which is an obstruction for representability by submanifolds.
As explained in~\cite[P.~771/2]{Hir2}, this technical formulation of Problem~8 was resolved by M.~Kreck, who performed a calculation
with the Atiyah--Hirzebruch spectral sequence. His result about integral Steenrod powers in particular implies that there is 
a manifold of dimension $10$ with a degree $7$ homology class not representable by a submanifold. However, Kreck's is 
purely an existence result, and no manifold is produced which has this property.

In~\cite{BHK}, C.~Bohr, B.~Hanke and I pointed out that there are explicit examples of $10$-manifolds, like the group manifold $Sp(2)$, which have 
degree $7$ homology classes not representable by submanifolds. This is detected by the reduced Steenrod powers in cohomology
with coefficients mod $3$, rather than with integral coefficients. The calculation of Steenrod powers required for $Sp(2)$ had been carried 
out at the beginning of the 1950s by A.~Borel and J-P.~Serre, so that the example of~\cite{BHK} was -- in some sense -- known when Hirzebruch
compiled the problem list~\cite{Hir1}, though clearly neither he nor Thom realized this. Of course the result of Borel and Serre was 
quite new at the time of~\cite{Hir1}, but was certainly well known in the mid 1980s when Kreck gave his solution recorded in~\cite{Hir2}.

\vfill\eject

\newpage

\section{Existence of $G$-structures and their integrability}\label{s:G}

Subsection~1.3 of~\cite{Hir1} discusses what are nowadays called $G$-structures on manifolds and their integrability.
In the terminology of $G$-structures $G$ denotes a subgroup of $GL(n,\R)$. However, in~\cite{Hir1} the general linear group 
itself is denoted by $G$, and subgroups are denoted by $H$, so that one has to speak about $H$-structures. We shall follow this 
terminology here, so as not to confuse readers when referring back  to~\cite{Hir1,Hir2}.

Problem~9, attributed to E.~Calabi, asks for a general criterion to determine whether a manifold $M^n$ admits an integrable 
$H$-structure, for some given closed subgroup $H\subset GL(n,\R)$. As explained in~\cite{Hir2}, many important cases
of this problem have been resolved, but the methods and techniques used are very different in the different cases, suggesting
that the general formulation attempted by Calabi may in fact not be appropriate.

I want to make a few comments here on some special cases, in part updating and expanding on Hirzebruch's commentary~\cite{Hir2}.

\medskip

\noindent
{\it (Almost) complex structures. \ }
For $n=2k$ and $H=GL(k,\C)$ an $H$-structure is an almost complex structure. By the Newlander--Nirenberg theorem,
the obstruction to integrability is precisely the Nijenhuis tensor. In real dimension $4$ there are many examples of 
almost complex manifolds without any complex structure. Such examples can be detected from many different points of view, e.g.
from the Enriques--Kodaira classification and Chern number inequalities, or using the fundamental group and the Albanese map,
or using gauge theory (Donaldson and Seiberg--Witten invariants).

In high dimensions, meaning real dimensions $\geq 6$, S.-T.~Yau has conjectured that every almost complex manifold should admit
a complex structure; see, for example, \cite[Problem~52]{Yau}. There is no real evidence in favor of Yau's speculation, but if it makes 
any sense at all, then it should perhaps be formulated in an even stronger form, as an $h$-principle: in dimensions $\geq 6$ every 
almost complex structure should be homotopic to an integrable one.

Note that not only are there almost complex four-manifolds without any complex structure, but there are also many complex surfaces 
that admit homotopy classes of almost complex structures that contain no integrable structure. 

\medskip

\noindent
{\it Symplectic structures. \ }
For $n=2k$ and $H=Sp(2k,\R)$ an $H$-structure is an almost symplectic structure, which is the same as an almost complex 
structure. One can think of an almost symplectic structure as a non-degenerate two-form, and the integrability condition is then
that this two-form be closed.

There are almost symplectic manifolds, like $S^1\times S^3$, or $S^6$, for which the structure of the cohomology ring obstructs the 
existence of a symplectic structure. Therefore, in order to have a chance to get a symplectic structure, one needs to assume not only 
that one has an almost symplectic manifold, but in addition it has to be cohomologically symplectic, or $c$-symplectic for short.
For $n=4$ this is still not enough for integrability. There are many examples of almost symplectic $c$-symplectic four-manifolds 
that cannot be symplectic, because of obstructions coming from Seiberg--Witten invariants, cf.~\cite{Bour}. In higher 
dimensions no further obstructions are known, but it is not known whether every almost symplectic structure on a $c$-symplectic 
manifold is homotopic to an integrable one. In any case, existence of symplectic structures is known in many cases; 
see the work of Gompf~\cite{Gompf} in particular.

\medskip

\noindent
{\it Pairs of complementary foliations. \ }
For $n=p+q$ and $H=GL(p,\R)\times GL(q,\R)$ an $H$-structure on $M$ is a splitting of the tangent bundle into a direct sum 
of subbundles or distributions of ranks $p$ and $q$ respectively. The question about the existence of such a splitting can often be answered 
in terms of algebraic topological invariants of $M$. Such an $H$-structure is integrable if and only if it is induced from a local 
product structure given by a pair of complementary foliations.

As mentioned in~\cite{Hir2}, if one considers the weaker requirement that only one of the two distributions be integrable, then a lot 
is known. This is the question of the existence of a foliation, of dimension $p$ say, on $M$, assuming that the tangent bundle of 
$M$ admits a rank $p$ subbundle. In many cases the integrability of all distributions (up to homotopy) has been proved by W.~P.~Thurston,
for example if $p=n-1$, see~\cite{ThurstonAnnals}, and also if $p=2$, see~\cite{ThurstonCMH}. 
(The case $p=1$ is trivial, of course.) In other cases there are additional obstructions 
coming from the so-called Bott vanishing theorem that forces the vanishing of certain characteristic classes of the normal bundle of a foliation.

Returning to the full integrability of $GL(p,\R)\times GL(q,\R)$-structures, it is still a very difficult problem to understand when a splitting 
of the tangent bundle can be induced by a pair of complementary foliations. Even in situations where both distributions are separately
homotopic to integrable ones, for example because of Thurston's results, it is very unclear whether they can simultaneously be made 
integrable in such a way that they remain complementary. It is certainly not possible to fix one foliation and then homotope the 
normal bundle to obtain a second, complementary, foliation. This problem appears already for $p=1$ and $q=2$, since there are 
circle bundles over surfaces which do not admit any horizontal foliation complementary to the fibers. Although the two-dimensional
horizontal foliation is homotopic to a foliation, that foliation will never be complementary to the fibers. Of course in this case one can 
just switch the r\^oles of the two distributions and argue that one makes the two-dimensional distribution integrable without worrying
about the complement since every one-dimensional distribution is integrable. This switching does not work already for $p=q=2$.
In this case all distributions are homotopic to integrable ones~\cite{ThurstonCMH}, but if we take for $M^4$ the non-trivial $S^2$-bundle over $S^2$,
then again there is no two-dimensional foliation complementary to the fibers of the fibration. If one just homotopes the horizontal
foliation to make it integrable (and no longer horizontal), then one does not know whether an integrable complement exists for 
the homotoped distribution. For general surface bundles over surfaces of positive genus the (non-)existence of a horizontal foliation
is an interesting problem that has attracted quite a bit of attention in recent years, but is still open.

\medskip

Next we come to a problem about $G$-structures defined by nonvanishing one-forms.
\begin{prob10}
A Pfaffian equation in $M^{2k+1}$ can be locally reduced to the normal form
$$
dx_{2k+1}+x_1dx_{k+1}+\ldots+x_rdx_{k+r}=0\ , \ \ \ \ 0\leq r\leq k\ ,
$$
where $r$ is the only local invariant, $2r+1$ being what is usually called the class of the local Pfaffian equation. The number $r$
being fixed, give a criterion to determine whether $M^{2k+1}$ admits a Pfaffian equation in the large whose class is $2r+1$ everywhere.
\end{prob10}

In more global terminology, the problem asks for criteria to determine whether, for given $r$, a manifold $M$ admits a 
one-form $\alpha$ for which $\alpha\wedge (d\alpha)^r$ is nowhere zero, but $\alpha\wedge (d\alpha)^{r+1}$ vanishes identically.
(Having a global defining form for the Pfaffian equation is no significant restriction, and can always be ensured by passage to 
a two-fold covering, if necessary.)
The problem also makes sense for even-dimensional manifolds, not just odd-dimensional ones. A slightly different problem would be 
to ask for $\alpha$ to be of constant class -- the notion of class for a one-form and for the Pfaffian equation it determines are 
slightly different. Chern's problem asks for the Pfaffian equation $\alpha=0$ to have constant class $2r+1$, which is equivalent 
to $\alpha$ having class $2r+1$ or $2r+2$, as a one-form, at every point of $M$. At points where $(d\alpha)^{r+1}=0$, the class 
of $\alpha$ is $2r+1$, and at points where $(d\alpha)^{r+1}\neq 0$, the class of $\alpha$ is $2r+2$.

For this problem, the cases corresponding to the two extreme values of $r$ have been studied the most. Firstly, if $r=0$, then 
one looks for nowhere zero one-forms $\alpha$ such that $\alpha\wedge d\alpha$ vanishes identically. By the Frobenius theorem,
this is equivalent to the existence of a codimension one foliation on $M$. It was proved by Thurston~\cite{ThurstonAnnals}
that the vanishing of the Euler characteristic, which is necessary for the existence of a codimension one distribution, is also 
sufficient for the existence of a foliation. In fact, in this case every distribution is homotopic to an integrable one.

Now the defining one-form for a codimension one foliation has class $1$ or $2$ at every point of $M$. A one-form whose class
is $1$ everywhere is nowhere vanishing and closed. The existence of such a form is equivalent to the existence of a foliation 
without holonomy, and, by D.~Tischler's theorem~\cite{T}, is also equivalent to the existence of a smooth fibration of $M$
over the circle. The existence of such a fibration is a much stronger condition than the vanishing of the Euler characteristic. For 
example, it implies the vanishing of the signature, and it imposes restrictions on the fundamental group.

Instead of the minimal value of $r$, let us now consider the maximal value. The case of even-dimensional manifolds is rather 
easier than that of odd-dimensional ones. So assume for the moment that $\dim (M)=2k+2$. Then the maximal $r$ for which
there can be a one-form $\alpha$ with $\alpha\wedge (d\alpha)^r$ nowhere zero, is $r=k$. For this value of $r$ we always 
have $\alpha\wedge (d\alpha)^{r+1}$ vanishing identically for dimension reasons. Now the kernel of such a one-form $\alpha$
is what is usually called a (cooriented) even contact structure. It is a result of D.~McDuff~\cite{McDuff} that these structures 
obey the $h$-principle, so that such structure exists for every reduction of the structure group of $T M$
from $GL(2k+2;\R)$ to $GL(k;\C)\times \{1\}\times\{ 1\}$.
In this case of even contact structures, the class of a defining one-form is everywhere $2k+1$ or $2k+2$. If $M$ is closed,
then the class cannot be $2k+2$ everywhere, since this would give an exact volume form on $M$. The class of $\alpha$
is $2k+1$ everywhere if and only if the holonomy of the characteristic foliation $C(\alpha)=\ker (\alpha\wedge (d\alpha)^k)$
is volume-preserving.

Returning to the case $\dim (M)=2k+1$ in the statement of the problem, the maximal value of $r$ for which $\alpha\wedge (d\alpha)^r$ 
can be nowhere zero is again $r=k$. In this case the kernel of $\alpha$ is a (cooriented) contact structure. Contact structures
are known to satisfy the $h$-principle only on open manifolds; see M.~Gromov~\cite{Gromov} and Y.~Eliashberg--N.~Mishachev~\cite{EM}.
Nevertheless, it is a classical result of R.~Lutz and J.~Martinet that on orientable three-manifolds every hyperplane distribution is 
in fact homotopic to a contact structure. Furthermore, a relationship between open book decompositions and contact structures
was discovered by W.~P.~Thurston and H.~E.~Winkelnkemper~\cite{TW}. E.~Giroux~\cite{Giroux} developed this into a precise 
correspondence between isotopy classes of contact structures and certain equivalence classes of open book decompositions.

On manifolds of higher dimension, such general results are not known.
A necessary condition for the existence of a contact structure is the existence of an almost contact structure, that is, a 
reduction of the structure  group of $T M$ from $GL(2k+1;\R)$ to $GL(k;\C)\times \{1\}$. Already in dimension $5$, this is a 
non-trivial constraint, but one that can be discussed effectively in terms of characteristic classes. More than twenty years ago,
H.~Geiges~\cite{Geiges1} proved that on simply connected $5$-manifolds every homotopy class of almost contact structures
does contain a contact structure. Although it was clear that the fundamental group does not obstruct the existence of contact 
structures, see~\cite{AcampoKot}, the generalization of Geiges's existence result took a very long time. After lots of partial 
results for special classes of manifolds, many of them obtained by Geiges and his coauthors, a complete existence result for 
arbitrary five-manifolds, and all homotopy classes of almost contact structures on them, was only announced very recently 
by R.~Casals, D.~Pancholi and F.~Presas~\cite{CPP} and by J.~B.~Etnyre~\cite{Etnyre}.

In dimensions $\geq 7$, there are few existence results, see however~\cite{Bourgeois,Geiges2,Geiges3}.
An extension of the correspondence between open book decompositions and contact structures, which is only known
in dimension three~\cite{TW,Giroux}, to higher dimensions was at one time announced by E.~Giroux and J.-P.~Mohsen.
The status of that, very technical, work is still unclear at the time of writing, and, therefore, some of the results in the literature 
that are based on it may have to be taken with a grain of salt. If that correspondence could indeed be put on a sound footing, then
it might form the basis for an inductive proof of a general existence theorem, where the induction hypothesis is 
applied on the spine of an open book decomposition, similarly to the application of the three-dimensional 
results in the proof of the existence theorem in dimension five in~\cite{Etnyre}.

\vfill\eject

\newpage

\section{Chern's embedding problems}\label{s:Chern}

Subsection~1.4 of~\cite{Hir1} states two embedding problems posed by S.-S.~Chern, the first one about smooth and the second one 
about isometric embeddings. One can think of both these problems as asking for optimal embedding dimensions 
for certain types of embeddings of manifolds in Euclidean spaces.

\begin{prob11}
Let $n+d(n)$ be the minimum dimension of an Euclidean space in which the $n$-dimensional real projective space $\PP^n$
can be differentiably imbedded. H.~Hopf proved that (even for topological imbedding) $d(n)>1$; Chern proved that $d(n)>2$ 
if $n\neq 2^k-1$ ($k\geq 2$) and $n\neq 2^k-2$ ($k\geq 2$). Can these bounds be improved? In particular, can $\PP^6$
be imbedded in an Euclidean space of $8$ dimensions?
\end{prob11}

As explained by Hirzebruch in his commentary~\cite[p.~774]{Hir2}, the explicit questions posed here were answered long ago:
yes, the bounds on embedding dimensions proved in the 1940s can often be improved, and no, $\PP^6$ does not embed in $\R^8$.
However, the more general question implicit in Chern's problem, namely to determine for {\it every} $n$ the optimal embedding
dimension for $\PP^n$, is still open. For example, it is now known that $\PP^6$ embeds in $\R^{11}$, but it is unknown whether
it embeds in $\R^{10}$.

This question has been a test case for many of the developments in algebraic topology and homotopy theory over the past 60 years. 
Unfortunately, the results that are known to date do not point to there being a clean answer that would be easy to state. 

In the early 1960s, M.~F.~Atiyah attempted a clean answer, by suggesting that the minimal embedding dimension for $\PP^n$ should be 
$2n-\alpha (n)+1$, where $\alpha (n)$ is the number of nonzero terms in the dyadic expansion of $n$, cf.~\cite{ES}. Minimality of this dimension 
was almost immediately disproved by M.~Mahowald~\cite{Maho}, with the first counterexample being $n=12$, where $2n-\alpha (n)+1=23$ 
and Mahowald gave an embedding in dimension $21$. However, it has turned out that the other half of Atiyah's suggestion was indeed 
correct: $\PP^n$ does embed in the Euclidean space of dimension $2n-\alpha (n)+1$. The first step in this direction was taken by D.~B.~A.~Epstein 
and R.~L.~E.~Schwarzenberger~\cite{ES}. Their work in particular gives the embedding of  $\PP^6$ in $\R^{11}$ mentioned above.
The general case, for $n$ odd only, was proved by B.~Steer~\cite{Steer}. Using Steer's result for the projective space of 
one more dimension does not quite resolve the case of even $n$, but the improvement (by one) needed for this was made 
by A.~J.~Berrick~\cite{Ber}.

Although dimension $2n-\alpha (n)+1$ always works for $\PP^n$, and is minimal for certain values of $n$, it is sometimes 
not minimal, as first shown by Mahowald~\cite{Maho}. However, there is a dearth of positive results, and most of the progress 
made beyond the work of Steer and Berrick has been negative, by proving non-embedding results for various smaller dimensions. 
Some of these results obtained before 1987 were summarized 
in the commentary~\cite[p.~774]{Hir2}. After 1987 many more incremental results have been obtained, with 
perhaps the most recent ones due to D.~M.~Davis in~\cite{Davis}. We refer the reader to this paper and to the literature quoted there,
as well as to the webpage~\cite{DavisTable} maintained by Davis that keeps track of the progress on this problem.

So far, no substitute for Atiyah's conjecture has emerged, meaning that a general solution is not in sight. Since people have concentrated
on specific, incremental, improvements for a long time, it is amazing that even for $\PP^6$ one does still not know whether the minimal 
embedding dimension is $9$, or $10$ or $11$, and no progress has been made on this particular case for about 50 years.

Finally, I cannot resist mentioning a far-reaching generalization of the correct part of Atiyah's conjecture: S.~Gitler~\cite{Git} conjectured
that {\it all} $n$-manifolds embed in $\R^{2n-\alpha (n)+1}$. As we saw above, this is now known for real projective spaces, but it is 
still open in general. 

\medskip

We now come to Chern's question about isometric embeddings.
\begin{prob12}
A closed surface $M$ in the $4$-dimensional Euclidean space $\EE^4$ has an induced Riemannian metric and hence a Gaussian 
curvature $K$. Does there exist such a surface for which $K<0$ everywhere? In particular, does there exist a surface of genus 
$2$ in $\EE^4$ for which $K<0$ everywhere?
\end{prob12}

In the 1987 update, Hirzebruch wrote~\cite[p.~775]{Hir2}: 

\begin{quote}
{\it Mir ist dar\"uber nichts bekannt.}
\end{quote}

The context for this problem was that, on the one hand, a closed surface in $\EE^3$ must have a point of positive Gaussian curvature
(by the maximum principle),
and, on the other hand, in the early 1950's the isometric embedding problem for Riemannian manifolds (in arbitrary codimension) was 
very much in the air. Shortly after Hirzebruch's problem list, J.~F.~Nash solved the embedding problem, and this implies of course that some 
$\EE^n$, with $n$ large, contains closed surfaces of negative Gaussian curvature, and even of constant negative Gaussian curvature. 
Eventually the isometric embedding dimension for Riemannian surfaces was reduced to $5$, see M.~Gromov~\cite[p.~298--303]{GromovBook}, 
showing in particular that $\EE^5$ contains surfaces of constant negative curvature for all genera $\geq 2$. In 1962, 
\`E.~R.~Rozendorn~\cite{Roz} constructed a closed surface of higher genus in $\EE^4$ with (non-constant) strictly negative Gaussian  
curvature, thus answering the first part of Problem~12. We refer to~\cite{Roz2} for a sketch of Rozendorn's construction. It is still an
open question whether one can reduce the genus to $2$, and whether one can make the curvature constant.

\vfill\eject

\newpage

\section{Dolbeault cohomology for almost complex manifolds}

Subsection~2.2 of~\cite{Hir1} discusses (almost) Hermitian metrics on closed almost complex manifolds. An almost complex structure 
on $M$ defines a decomposition $T M\otimes_{\R}\C=T^{1,0}M\oplus T^{0,1}M$, which in turn induces a decomposition of 
complex differential forms into $(p,q)$-types, just as in the complex case. One then defines differential operators $\del$, respectively $\bar\del$,
as the components of the exterior derivative $d$ that raise $p$ respectively $q$ by one. On functions or zero-forms one always has 
$d=\del+\bar\del$, but this is no longer true on forms of higher degree as there are additional summands arising from the Nijenhuis 
tensor. Nevertheless, one can consider the $L^2$-adjoints of these operators with respect to an Hermitian metric and define the
$\bar\del$-Laplacian
$$
\square=\bar\del^*\bar\del+\bar\del \bar\del^* \ .
$$
In~\cite[Subsection~2.2]{Hir1} the space $H^{p,q}$ is defined to be the kernel of $\square$ in the space of $(p,q)$-forms, and $h^{p,q}$ 
is its dimension.

\begin{prob20}
Let $M^n$ be an almost-complex manifold. Choose an Hermitian structure and consider the numbers $h^{p,q}$ defined as above.
Is $h^{p,q}$ independent of the choice of the Hermitian structure? If not, give some other definition of the $h^{p,q}$ of $M^n$ which 
depends only on the almost-complex structure and which generalizes the $h^{p,q}$ of a complex manifold.
\end{prob20}

On a complex manifold the $H^{p,q}$ are the Dolbeault cohomology groups, isomorphic to the sheaf cohomology of the bundle of 
holomorphic $p$-forms. Therefore, in that case the $h^{p,q}$ are indeed independent of the metric.

As pointed out in the commentary in~\cite[p.~779]{Hir2}, the Atiyah-Singer index theorem implies that for fixed $p$ the index of 
$\bar\del+\bar\del^*$ is the $p$-component of the Hirzebruch--Todd genus, and is therefore metric-independent without assuming
integrability. However, the operator $\square$ is different from $(\bar\del+\bar\del^*)^2$ in the non-integrable case.

There seems to have been no progress at all on this problem, which asks for a development of harmonic Dolbeault theory 
on arbitrary almost complex manifolds. Such a theory could be very useful, particularly for geometrically interesting
almost complex structures, like those tamed by a symplectic form. In this context an attempt to develop harmonic theory 
was made by S.~K.~Donaldson in~\cite{DonDur}.

M.~Verbitsky~\cite{Ver} has developed Hodge theory on strictly nearly K\"ahler $6$-manifolds. These are certain special 
almost complex manifolds with a Hermitian metric with respect to which the covariant derivative of the almost
complex structure is skew-symmetric. For these special Hermitian metrics Verbitsky obtained a Hodge decomposition of 
the harmonic forms for the usual Riemannian Laplacian $\Delta=dd^*+d^*d$ into harmonic forms of pure $(p,q)$-types.
However, $\Delta$ is different from the $\bar\del$-Laplacian $\square$ considered in Problem~20.

\vfill\eject

\newpage

\section{Existence and classification of complex structures}

Three problems that are still open concern complex structures on certain simple manifolds. We first state the two problems concerning 
complex projective spaces\footnote{In Section~\ref{s:Chern} $\PP^n$ denoted real projective $n$-space, in this section $\PP_n$ denotes complex
projective $n$-space.}.

\begin{prob13}
Does there exist a complex structure on $\PP_3$ with vanishing second Chern class? (Such a complex structure cannot carry a 
K\"ahlerian metric, see [\ldots] Problem 28*.) 
\end{prob13}

Hirzebruch pointed out that an almost complex structure on $\PP_3$ with $c_2=0$ does exist, and is obtained by blowing up a point in
$S^6$ equipped with an almost complex structure. A negative answer to Problem~13 would imply that $S^6$ does not have a complex
structure. This problem is still considered to be wide open, although solutions have been claimed repeatedly over the years.

\medskip

\begin{prob28}
Consider the complex projective space $\PP_n$ as a differentiable manifold with the usual differentiable structure. Determine all 
complex structures of $\PP_n$. {\rm (See Problem~14.)} In particular (28*), determine all complex structures of $\PP_n$ 
which can carry a K\"ahler metric.
\end{prob28}

Problem~14 that is alluded to here asked about all the possible Chern classes of almost complex structures on $\PP_n$.
As described in~\cite[p.~776]{Hir2}, that problem was essentially solved by M.~Puschnigg~\cite{Puschnigg}.

Hirzebruch also explained in~\cite[p.~782]{Hir2} how Problem 28* is solved. If one has a K\"ahlerian complex structure on $\PP_n$,
then by the Kodaira embedding theorem it is projective algebraic. A combination of results of Hirzebruch and Kodaira proved 
in the mid-1950s with the later result of Yau ruling out structures with ample canonical bundle shows that the standard complex
structure is the only algebraic one.

Both Problem~13  and Problem~28 are open for non-K\"ahlerian complex structures as soon as $n\geq 3$. For $n=2$ one 
does know that any complex structure on $\PP_2$ is automatically K\"ahler; see~\cite{Buch,Lam}. As explained above,
it is then standard. 

If Yau's conjecture mentioned in Section~\ref{s:G} above were true, then $S^6$ would be a complex manifold, and $\PP_3$ 
would admit at least one non-K\"ahlerian complex structure. If the strong version of the conjecture were true, in the form 
of an $h$-principle, then by the classification of almost complex structures due to Puschnigg mentioned above, $\PP_n$ 
would have infinitely many non-homotopic and therefore non-deformation-equivalent complex structures for every $n\neq 1, 2, 4$. 

\medskip

This brings us to a specifically four-dimensional problem:
\begin{prob25}
Determine all complex structures of $S^2\times S^2$.
\end{prob25}

After the complex projective spaces, this is the next natural example to look at. In his thesis, Hirzebruch had considered the complex
structures on $S^2\times S^2$ that we now call the (even) Hirzebruch surfaces\footnote{This was the subject of a lecture by Hirzebruch at
the Cornell conference, cf.~\cite{Steenrod}.}. They show that $S^2\times S^2$ has infinitely many
non-biholomorphic complex K\"ahler structures. These structures are all ruled and projective algebraic. Problem~25 asks whether there are
any other complex structures on $S^2\times S^2$ beyond the Hirzebruch surfaces. 
If we fix the standard smooth structure on $S^2\times S^2$, then the answer is that there are no other complex structures on it.
If we allow all manifolds homeomorphic to $S^2\times S^2$, then the answer is still not known, and in this sense the problem is 
still open.

To explain all this, 
let us consider more generally all compact complex surfaces $X$ with the Betti numbers and the intersection form of $S^2\times S^2$,
without assuming for the moment that $X$ is simply connected, and with no assumption on the underlying smooth structure.
Then $X$ has to be K\"ahlerian. This can be deduced from classification results of K.~Kodaira, but a direct argument is now 
available, due to N.~Buchdahl~\cite{Buch} and A.~Lamari~\cite{Lam}. It then follows from the Kodaira embedding theorem that $X$ 
is projective algebraic. Its Chern numbers are $c_1^2=8$ and $c_2=4$. Since the intersection form is even, $X$ is minimal.
Now the Enriques-Kodaira classification implies that $X$ is either ruled, in which case it is a Hirzebruch surface, or is a minimal 
surface of general type.

At the time when Hirzebruch wrote the commentary in~\cite{Hir2}, it was an important problem to decide whether a surface of general type 
could ever be diffeomorphic to a rational surface. At that time, the so-called Barlow surface was the only known example of a simply 
connected minimal surface of general type with geometric genus $p_g=0$, equivalently with $b_2^+=1$.
In fact, by M.~H.~Freedman's classification up to homeomorphism, the Barlow surface was known to be homeomorphic to the 
$8$-fold blowup of the complex projective plane. In 1988 I proved that the Barlow surface is not diffeomorphic to a 
rational surface~\cite{Invent}. There were then several attempts to generalize this result to other surfaces of general type, although
no other ones were known to exist that were even homeomorphic to rational surfaces. Finally, R.~Friedman and Z.~Qin~\cite{FQ}
proved that no surface of general type can be diffeomorphic to a rational surface, in particular not to $S^2\times S^2$.
(This was an important result in Donaldson theory, but became fairly straightforward after the advent of 
Seiberg--Witten theory~\cite{FM,Wi}.) Thus, the only complex surfaces {\it diffeomorphic} to $S^2\times S^2$ are the even Hirzebruch surfaces.

Starting in 2006, J.~Park, together with several coauthors, constructed new examples of surfaces of general type that are homeomorphic 
to the complex projective plane blown up in $5$, $6$ or $7$ points. We refer to Park's survey in~\cite{Park}. By~\cite{FQ,FM,Wi} 
these surfaces are not diffeomorphic to rational surfaces. No smaller simply
connected examples have been found so far, in particular none that have the Betti numbers of $S^2\times S^2$. 

If we do not insist on simple connectivity and consider arbitrary surfaces $X$ of general type with $p_g=0$, $c_1^2=8$ and $c_2=4$, then examples
have been known for a long time, and were mentioned already by Hirzebruch in the commentary in~\cite{Hir2}. However, all known
examples are uniformized by the polydisk $\HH^2\times\HH^2$, and therefore have infinite fundamental group. Some authors have 
speculated about a general uniformization result, that would force all surfaces with these numerical invariants to be quotients of the polydisk.
If that could be proved, then any compact complex surface {\it homeomorphic} to $S^2\times S^2$ would be a Hirzebruch surface.

The speculations about a uniformization result are motivated by the parallel case of ball quotients. 
Through the work of Y.~Miyaoka and S.-T.~Yau one knows that any surface with $c_1^2=3c_2>0$ is either the projective plane or a 
compact quotient of its non-compact dual $\C H^2$. For the case $c_1^2=2c_2>0$, equivalently zero signature and positive Euler 
characteristic, some partial uniformization results have been obtained using gauge theory. First, using Yang--Mills theory, 
C.~T.~Simpson~\cite{Simp} proved that any K\"ahler surface with $c_1^2=2c_2>0$ having the additional 
property that the tangent bundle is the direct sum of two line bundles of negative degrees is uniformized by the 
polydisk. Second, I proved that a surface $X$ of general type with $c_1^2=2c_2>0$ for which the Seiberg--Witten
invariants of the underlying smooth manifold do not vanish for the orientation opposite to that induced by the complex 
structure is also uniformized by the polydisk~\cite{BLMS}. Thus, if the underlying manifold
has an orientation-reversing self-diffeomorphism, or if it is complex or symplectic for the opposite orientation
in some other way, then $X$ is uniformized by the polydisk and thus has infinite fundamental group. In particular 
it cannot be homeomorphic to $S^2\times S^2$.

That the additional assumptions in these uniformization results cannot be dispensed with follows from the existence 
of several important examples that disprove an unconditional statement for surfaces of general type with zero signature.
Among surfaces of general type there are simply connected ones satisfying $c_1^2=2c_2>0$. They can even be chosen to 
be spin by a result of B.~Moishezon and M.~Teicher~\cite{MT}.
These examples are homeomorphic to connected sums of very many copies of $S^2\times S^2$.
Very recently some new surfaces of general type with zero signature have been found that
are not uniformized by the polydisk~\cite{CMR}, although they have infinite fundamental group. In this case $c_1^2=2c_2=16$.

\vfill\eject

\newpage

\section{Topological invariance of characteristic numbers of algebraic varieties}

In the following problem $h^{p,q}$ denotes the Hodge numbers.

\begin{prob31}
Are the $h^{p,q}$ and the Chern characteristic numbers of an algebraic variety $V_n$ topological invariants of $V_n$? If not,
determine all those linear combinations of the $h^{p,q}$ and the Chern characteristic numbers which are topological invariants.
\end{prob31}

This problem was motivated by the Hirzebruch--Riemann--Roch theorem, which at the time of~\cite{Hir1} was only known for 
complex projective varieties. Of course the problem is interesting more generally, especially for K\"ahler manifolds, where the 
Hodge numbers have the same properties and symmetries as on algebraic varieties.

It is clear from the text in~\cite{Hir1} that Hirzebruch was well aware of the fact that Chern numbers of almost complex structures 
would probably not be topological invariants of the underlying manifold, since there are too many almost complex structures with 
different Chern classes. The term ``topological invariant'' can be interpreted in several different ways, for example one can take it 
to mean homeomorphism- or diffeomorphism-invariance, and one can fix the orientation, or not. It is obvious that the  top Chern 
number is a topological invariant in every sense of the word, since it is the Euler number, and Hirzebruch~\cite{Hir1} knew that 
the Pontryagin numbers are oriented diffeomorphism invariants.
In fact, by S.~Novikov's result, they are oriented homeomorphism invariants. By the Hodge decomposition of the 
cohomology, there are certain linear combinations of Hodge numbers that reduce to Betti numbers, and are therefore 
topological invariants.

Up until~\cite{Hir2}, the only further information in the direction of Problem~31 was the fact that $c_1^5$  is not a diffeomorphism 
invariant of complex projective $5$-folds. An example showing this had been found by Borel and Hirzebruch in 1958. Of course 
one can take products of this example with other varieties to generate examples of diffeomorphic varieties in higher dimensions 
which have distinct Chern numbers. These examples say nothing about the Hodge numbers, and, as recently as 2010,
algebraic geometers were asking each other in internet discussion groups whether the Hodge numbers of K\"ahler manifolds 
are diffeomorphism invariants.

After  some earlier, partial, results in small dimensions, see~\cite{JTop}, I completely solved Problem~31 for Chern numbers in 2009, 
see~\cite{PNAS,Adv}. I proved that a rational linear combination of Chern numbers is an oriented diffeomorphism invariant of smooth 
complex projective varieties if and only if it is a linear combination of the Euler and Pontryagin numbers. If one does not fix the orientation, 
then in complex dimension $n\geq 3$ a rational linear combination of Chern numbers is a diffeomorphism invariant of smooth complex 
projective varieties if and only if it is a multiple of the Euler number $c_n$. In complex dimension $2$ both Chern numbers $c_2$ and 
$c_1^2$ are diffeomorphism-invariants of complex-algebraic surfaces, see~\cite{BLMSorient,JTop}. Except in complex dimension $2$, 
the results are the same if one considers homeomorphism invariants instead of diffeomorphism invariants.

To prove these results one needs a supply of diffeomorphic projective algebraic varieties with distinct Chern numbers.
Using the unitary bordism ring for bookkeeping, one can reduce the problem to the construction of certain special basis 
sequences for this bordism ring tensored with $\Q$. Suitable basis sequences involving formal differences of diffeomorphic 
varieties are constructed in~\cite{Adv} by considering the pairs of orientation-reversingly homeomorphic algebraic surfaces 
obtained in~\cite{MAorient} and 
algebraic projective space bundles over them, so that the total spaces become diffeomorphic, but have distinct Thom--Milnor 
numbers (which are special combinations of Chern numbers).

The full solution to Problem~31 was not obtained in~\cite{PNAS,Adv}, because I did not discuss Hodge numbers systematically there.
I did understand that the only linear relations between Hodge and Chern numbers are the Hirzebruch--Riemann--Roch 
relations $\chi(X;\Omega^p)=\Td_p(X)$, where the left-hand side is a combination of Hodge numbers and the right-hand 
side is a combination of Chern numbers. While the examples in~\cite{PNAS,Adv} do show that certain Hodge numbers 
are not diffeomorphism-invariant, there was no way to do a good bookkeeping of all possible linear combinations because 
the Hodge numbers are not bordism invariants.

These issues were resolved in recent joint work with S.~Schreieder~\cite{KS}, in which we studied the Hodge ring of K\"ahler manifolds.
This is a ring that keeps track of Hodge numbers in the same way that the unitary bordism ring keeps track of Chern numbers.
Refining the results of~\cite{PNAS,Adv} on the completeness of the Hirzebruch--Riemann--Roch relations we also 
studied a mixed Chern--Hodge ring, in which two equidimensional K\"ahler manifolds have the same image if and only if 
the have the same Hodge and Chern numbers. Using the Chern--Hodge ring, we completely solved Problem~31
for mixed linear combinations of Hodge and Chern numbers. The answer is, cf.~\cite[Theorem~4]{KS}:

\begin{quote}
{\it A rational linear combination of Hodge and Chern numbers of smooth complex projective varieties is
\begin{itemize}
\item[(1)] an oriented homeomorphism or diffeomorphism invariant if and only if it reduces to a linear combination of the Betti and Pontryagin 
numbers after perhaps adding a suitable combination of the $\chi (\Omega^p)-\Td_p$, and
\item[(2)] an unoriented homeomorphism invariant in any dimension, or an unoriented diffeomorphism invariant in dimension $n\neq 2$,
if and only if it reduces to a linear combination of the Betti numbers after perhaps adding a suitable combination of the $\chi(\Omega^p)-\Td_p$.
\end{itemize}
}
\end{quote}

\bigskip

\end{document}